\documentclass[10pt,twoside]{amsart}
\usepackage{amsmath}
\usepackage{amssymb}

\usepackage{graphicx}

\usepackage[colorlinks=true, urlcolor=blue,bookmarks=true,bookmarksopen=true, citecolor=blue]{hyperref}

\numberwithin{equation}{section}

\newtheorem{theorem}{Theorem}[section]
\newtheorem{definition}[theorem]{Definition}

\newtheorem{corollary}[theorem]{Corollary}

\newtheorem{lemma}[theorem]{Lemma}
\newtheorem{proposition}[theorem]{Proposition}

\usepackage{mathrsfs}         
\usepackage[all,pdf]{xy}
\usepackage{tikz}
\usetikzlibrary{arrows}

\newtheorem{question}[theorem]{Question}

\def \bb{\mathbb}
\def \mb{\mathbf}
\def \mc{\mathcal}

\def \ms{\mathscr}

\renewcommand\tilde{\widetilde}

\def \CC{{\bb{C}}}
\def \CP{{\bb{CP}}}

\def \HH{{\bb H}}
\def \JJ{{\bb{J}}}

\def \RR{{\bb{R}}}
\def \RP{{\bb{RP}}}

\def \ZZ{{\bb{Z}}}

\def \MMC{{\mc M}}

\def \PPC{{\mc P}}

\def \({\left(}
\def \){\right)}
\def \<{\langle}
\def \>{\rangle}

\def \bar{\overline}

\def \dsum{\oplus}

\def \inter{\cap}

\def \union{\cup}
\def \vargeq{\geqslant}
\def \varleq{\leqslant}

\def \xto{\xrightarrow}

\def \Ham{{\rm Ham}}

\def \Homeo{{\rm Homeo}}

\def \img{{\rm img }}

\renewcommand{\1}{1\!\!\!1}

\def \qed{\hfill $\square$ \vspace{0.03in}}

\begin{document}

\title{An example concerning Hamiltonian groups of self product, II}

\author{Shengda Hu}
\address{Department of Mathematics, Wilfrid Laurier University, 75 University Ave. West, Waterloo, Canada}
\email{shu@wlu.ca}

\author{Fran\c{c}ois Lalonde}
\address{D\'epartement de math\'ematiques et de Statistique, Universit\'e de Montr\'eal, C.P. 6128, Succ. Centre-ville, Montr\'eal H3C 3J7, Qu\'ebec, Canada}
\email{lalonde@dms.umontreal.ca}

\abstract
We describe the natural identification of $FH_*(X \times X, \triangle; \omega \dsum -\omega)$ with $FH_*(X, \omega)$. Under this identification, we show that the extra elements in $\Ham(X \times X,  \omega \dsum -\omega)$ found in \cite{HuLalonde3}, for $X = (S^2 \times S^2, \omega_0 \dsum \lambda \omega_0)$ for $\lambda > 1$, do not define new invertible elements in $FH_*(X, \omega)$.
\endabstract

\maketitle


\vspace{.08in} \noindent \textbf{Keywords}: Lagrangian submanifolds, Hamiltonian group, Seidel elements.

\section{Introduction}

Let $M$ be a symplectic manifold with an anti-symplectic involution $c$, such that $L$ is the Lagrangian submanifold fixed by $c$. 
For any map $u : (\Sigma, \partial \Sigma) \to (M, \triangle)$, where $\Sigma$ is a manifold with boundary, we define $v: \Sigma \union_{\partial} \bar \Sigma \to X$ by
$$v|_{\Sigma} = p_1 \circ u \text{ and } v|_{\bar\Sigma} = p_2 \circ u,$$
where $\bar \Sigma$ is $\Sigma$ with the opposite orientation. For any map $v : \Sigma \union_{\partial} \bar \Sigma \to X$ we obtain the corresponding map $u : (\Sigma, \partial \Sigma) \to (M, \triangle)$ by
$$u(x) = (v(x), v(\bar x)),$$
where $\bar x$ denotes $x \in \bar \Sigma$.
We use $\delta$ to denote the map $v \mapsto u$. 

For $M = X \times X$ and involution switching the factors, then $L = \triangle$. Let $\delta_k := p_k\circ \delta$ where $p_k$ is the projection to the $k$-factor, then it induces a map of Floer homologies. We show in \S\ref{Floer}
\begin{lemma}\label{lem:FHrelation}
There is a commutative diagram of isomorphisms of the Floer homologies:
\begin{equation*}
\text{
\xymatrix{
&& FH_*(M, \triangle) \ar[lld]^{\simeq}_{\delta_1} \ar[rrd]^{\delta_2}_{\simeq} && \\
FH_*(X) \ar[rrrr]^{\tau}_{\simeq} &&&& FH_*(\bar X)
}
}
\end{equation*}
\end{lemma}

In \cite{HuLalonde}, with proper assumptions, we described a construction of Lagrangian Seidel element from a path of Hamiltonian diffeomorphisms. In particular, a loop $\gamma$ in $\Ham(M)$ defines a Lagrangian Seidel element $\Psi^L_{\gamma} \in FH_*(M,L)$, where $L$ is a Lagrangian submanifold. The Albers' map 
$$\ms A : FH_*(M) \to FH_*(M,L)$$
whenever well-defined, for example when $L$ is monotone, relates the Seidel elements $\Psi^M_\gamma \in FH_*(M)$ to $\Psi^L_\gamma$. Let $S_M$ denote the image of Seidel map $\Psi^M : \pi_1\Ham(M) \to FH_*(M)$ and $S_L$ that of $\Psi^L : \pi_1(\Ham(M), \Ham_L(M)) \to FH_*(M, L)$, where $\Ham_L(M)$ is the group of Hamiltonian diffeomorphisms preserving $L$ which restrict to isotopies on $L$, then 
$$\ms A(S_M) \subseteq S_L$$
\begin{question}\label{quest:seidelinclusion}
When all the terms involved is well defined, is the inclusion $\ms A(S_M) \subseteq S_L$ (in general) proper?
\end{question}
\noindent
An affirmative answer to this question would imply an affirmative answer to the open question about the non-triviality of $\pi_0\Ham_L(M)$.

For the case of $L = \triangle$, 
since $\delta_1$ is an isomorphism, the inclusion is equivalent to $\delta_1\ms A(S_M) \subseteq \delta_1(S_\triangle)$ as subsets of $FH_*(X)$. 
In \S \ref{elements}, we show that $S_X \subseteq \delta_1\ms A(S_M)$. More precisely,
\begin{theorem}\label{thm:splitSeidel}[Corollary \ref{coro:splitSeidel}]
Let $\gamma \in \pi_1\Ham(X)$. It naturally lifts to a \emph{split} element $\gamma_+
\in \pi_1\Ham(M)$, and we have $\delta_1 \ms A(\Psi^M_{\gamma_+}) = \Psi^X_\gamma$.
\end{theorem}
\noindent
As a corollary, it shows that the natural map $\pi_1\Ham(X) \times \pi_1\Ham(\bar X) \to \pi_1\Ham(M)$ is injective. In light of this result, we pose the following question, which is related to Question \ref{quest:seidelinclusion} for the special case of diagonal.
\begin{question}\label{quest:diagonalinclusion}
Is any inclusion in the sequence $S_X \subseteq \delta_1\ms A(S_M) \subseteq \delta_1(S_\triangle)$ proper? 
\end{question}
\noindent
For $X = S^2 \times S^2$ as in \cite{HuLalonde3}, we show in \S \ref{elements} that the image under $\delta_1\ms A$ of the extra Seidel elements found in \cite{HuLalonde3} is contained in $S_X$. 

\vspace{0.1in}\noindent
{\bf Acknowledgement.} S. Hu is partially supported by an NSERC Discovery Grant.

\section{Identification of Floer homologies}\label{Floer}
\subsection{Notations}Let 
$$D^2_+ = \{z\in \CC : |z| \varleq 1, \Im z \vargeq 0\},$$
$\partial_+$ denote the part of boundary of $D^2_+$ on the unit circle, parametrized by $t \in [0,1]$ as $e^{i\pi t}$, and $\partial_0$ the part on the real line, parametrized by $t \in [0,1]$ as $2t -1$.

Let $(M,L)$ be a pair of symplectic manifold and a Lagrangian submanifold, and $\Omega$ is the symplectic form. For $\beta \in \pi_2(M,L)$, $\mu_L(\beta)$ denotes its Maslov number, and $\Omega(\beta)$ its symplectic area. The space of paths in $M$ connecting points of $L$ is
$$\PPC_L M = \{l : ([0,1],\partial[0,1]) \to (M, L), [l] = 0 \in \pi_1(M, L)\}$$ and 
the corresponding covering space with covering group $\Gamma_L = \pi_2(M,L)/ (\ker \omega \inter \ker \mu_L)$ is
$$\widetilde \PPC_L M = \{[l, w] : w : (D^2_+; \partial_+, \partial_0) \to (M; l,L)\}$$
where $(l, w) \sim (l',w') \iff l=l'$ and $\omega(w\#(-w')) = \mu_L(w\#(-w'))$. The space of contractible loops in $M$ parametrized by $\RR/\ZZ$ is denoted $\Omega(M)$ and the corresponding covering space with covering group $\Gamma_\omega = \pi_2(M)/(\ker \omega \inter \ker c_1)$ is given by
$$\widetilde \Omega(M) = \{[\gamma, v] : v : (D^2, \partial D^2) \to (M, \gamma)\}$$
where $(\gamma, v) \sim (\gamma', v') \iff \gamma = \gamma'$ and $\omega(v\#(-v')) = c_1(v\#(-v'))$. Here, $\partial D^2$ is parametrized as the unit circle in $\CC$ by $\{e^{2\pi i t}: t \in [0,1]\}$, and $c_1=c_1(TM)$ in some compatible almost complex structure. We denote the space of loops in $M$ parametrized by $\RR/T\ZZ$ and the corresponding covering space as $\Omega^{(T)}(M)$ and $\widetilde \Omega^{(T)}(M)$ respectively, thus $\Omega(M) = \Omega^{(1)}(M)$ and $\widetilde \Omega(M) = \widetilde \Omega^{(1)}(M)$.

Let $H : [0,1] \times M \to \RR$ be a time-dependent Hamiltonian function, which defines on $\tilde \PPC_L M$ the action functional
$$a_{H}([l, w]) = -\int_{D^2_+}w^*\omega + \int_{[0,1]} H_t(l(t)) dt,$$
where we use the convention $dH = - \iota_{X_H} \omega$ for the Hamiltonian vector fields. Similarly, a time dependent Hamiltonian function $K$ for $t \in \RR/T\ZZ$ defines an action functional $a_{K}$ on $\tilde \Omega^{(T)}(M)$. We will not distinguish notations for the two types of action functionals when it is clear from the context which one is under discussion.

Given the time dependent Hamiltonian function $H$, let $\tilde l \in \tilde \PPC_L M$ such that $l$ is a connecting orbit for $H$, then $\mu_H(\tilde l)$ denotes the corresponding Conley-Zehnder index. Similarly, for the time dependent Hamiltonian functino $K$, let $\tilde \gamma \in \widetilde \Omega^{(T)}(M)$ such that $\gamma$ is a periodic orbit for $K$, then $\mu_K(\tilde \gamma)$ denotes the corresponding Conley-Zehnder index. The following relations hold
$$\mu_H(\tilde l) - \mu_H(\tilde l') = \mu_L(w \#(-w')) \text{ and } \mu_K(\tilde \gamma) - \mu_K(\tilde \gamma') = c_1(u \#(-v'))$$
where $l = l'$ and $\gamma = \gamma'$.

\subsection{Doubling construction}
First we describe the doubling construction when the Lagrangian submanifold is the fixed submanifold of an anti-symplectic involution. It applies in this case since the diagonal $\triangle$ is the fixed submanifold of the involution of switching the two factors. 

Let $c : M \to M$ be an anti-symplectic involution and $L \subset M$ be the fixed submanifold of $\tau$, then it is a Lagrangian submanifold. We'll use $(\HH, \JJ)$ to denote a pair of $2$-periodical Hamiltonian functions and compatible almost complex structures, i.e.
$$\HH : \RR / 2 \ZZ \times M \to \RR \text{ and } \JJ = \{\JJ_t\}_{t \in \RR / 2\ZZ}.$$
\begin{definition}\label{involution:symdata}
The pair $(\HH, \JJ)$ is \emph{$c$-symmetric} if it satisfies
$$\HH_t(x) = \HH_{2-t}(c(x)) \text{ and } \JJ_t(x) = -dc\circ \JJ_{2-t}\circ dc.$$
For such a pair, we define the halves $(H, \mb J) := (\HH_t, \JJ_t)_{t\in [0,1]}$ and 
$$(H', \mb J') : = (\HH_{1-t}\circ c, -dc \circ \JJ_{1-t} \circ dc)_{t \in [0, 1]} = (\HH_{t+1}, \JJ_{t+1})_{t \in [0, 1]}.$$
\end{definition}

The doubling map $\delta$ described in the introduction is a special case of the following construction for a symplectic manifold with an anti-symplectic involution:
\begin{definition}\label{involution:doubledmap}Let $u: (\Sigma, \partial \Sigma) \to (M, L)$ be a map from a manifold $\Sigma$ with boundary $\partial \Sigma$, the \emph{doubled map} is given by:
$$v : \Sigma \union_\partial \bar\Sigma \to M : v|_{\Sigma} = u \text{ and } v|_{\bar\Sigma} = c \circ u,$$
where $\bar\Sigma$ is $\Sigma$ with the opposite orientation. We also write $\delta(u) := v$ which gives the \emph{doubling map} between the spaces of continuous maps:
$$\delta: Map(\Sigma, \partial \Sigma; M, L) \to Map(\Sigma\union_\partial \bar\Sigma; M).$$
\end{definition}
\noindent
In particular, we have the map between the space of paths in $(M, L)$ and loops of period $2$ in $M$, 
as well as their covering spaces:
$$\delta : \PPC_L M \to \Omega^{(2)}(M) \text{ and } \delta: \widetilde \PPC_L M \to \widetilde \Omega^{(2)}(M)$$

Let $(\HH, \JJ)$ be a $c$-symmetric pair
and $\{\phi_t\}_{t \in [0,2]}$ the Hamiltonian isotopy generated by $\HH$, then
\begin{equation}\label{involution:isodouble}
\phi_t = c \circ \phi_{2-t} \circ \phi_2^{-1} \circ c \Longrightarrow (c\circ \phi_2)^2 = \1.
\end{equation}
Let $(H, \mb J)$ and $(H', \mb J')$ be the two halves of $\HH$, then
$$H_t = H'_{1-t}\circ c \text{ and } J_t = -dc\circ J'_{1-t} \circ dc,$$
Let $\phi'_t$ denote the Hamiltonian isotopy generated by $H'$, then
$$\phi'_t = c \circ \phi_{1-t} \circ \phi_1^{-1} \circ c$$ It follows that if $l$ is a Hamiltonian path generated by $H$ connecting $x, y \in L$, then $l'(t) := c\circ l(1-t)$ is a Hamiltonian path generated by $H'$ connecting $y, x \in L$, and the double $\gamma = \delta(l)$ is a periodic orbit for $\HH$. This correspondence lifts to the covering spaces and the following holds.
\begin{lemma}\label{involution:actionfunct}
For $\tilde l \in \widetilde \PPC_LM$ let $\tilde \gamma = \delta(\tilde l)$, then 
$$a_\HH(\tilde \gamma) = 2a_H(\tilde l) = 2 a_{H'}(\tilde l').$$
Moreover, if $\tilde l$ is a critical point of $a_H$ then $\tilde \gamma$ is a critical point of $a_\HH$. If $\tilde \gamma$ is non-degenerate, then $\tilde l$ is as well. A Floer trajectory for $a_H$ is taken to a Floer trajectory for $a_\HH$ by $\delta$, which converges to the corresponding critical points when the trajectory has finite energy.
\qed
\end{lemma}

A result from \cite{HuLalonde} relates the Conley-Zehnder indices of connecting paths generated by $H$ and $H'$.
\begin{lemma}[Lemma 5.2 of \cite{HuLalonde}]\label{lem:reflectionindex}
Let $\tilde l$ and $\tilde l'$ be respective critical points of $a_H$ and $a_{H'}$ as above. Then $\mu_H(\tilde l) = \mu_{H'}(\tilde l')$. \qed
\end{lemma}

\subsection{Index comparison}
We briefly recall the definition of Conley-Zehnder index using the Maslov index of paths of Lagrangian subspaces as in Robbin-Salamon \cite{RobbinSalamon2}. Let $\tilde l = [l, w]$ be a non-degenerate critical point of $a_H$. Then $w : D^2_+ \to M$ and $l = \partial w$ is a Hamiltonian path. There is a symplectic trivialization $\Phi$ of $w^*TM$ given by
$\Phi_z : T_{w(z)} M \to \CC^n$ with the standard symplectic structure $\omega_0$ on $\CC^n$. Furthermore, we require that $\Phi_r(T_{w(r)}L) = \RR^n$, for $r \in [-1, 1] \subset D^2_+$. Then the linearized Hamiltonian flow $d\phi_t$ along $l$ defines a path of symplectic matrices
\begin{equation}\label{eq:symplecticpath}
E_t = \Phi_{e^{i\pi t}} \circ d\phi_t \circ \Phi_1^{-1} \in Sp(\CC^n) 
\end{equation}
Then the Conley-Zehnder index of $\tilde l$ is given by
$$\mu_H(\tilde l) = \mu (E_t \RR^n, \RR^n)$$
where $\mu$ is the Maslov of paths of Lagrangian subspaces introduced in \cite{RobbinSalamon2}.

We continue with the notations of Lemma \ref{involution:actionfunct}. 
\begin{proposition}\label{invol:indexdiff} 
Suppose that all the critical points involved are non-degenerate,
then
\begin{equation}\label{involution:indexdiffeq} \mu_{H}(\tilde l) + \mu_{H'}(\tilde l') - \mu_{\HH}(\tilde\gamma) = \frac{1}{2} sign(Q),\end{equation}
where $Q(\bullet, *) = \Omega((\1 - d\phi_2)\bullet, dc(*))$ is a quadratic form on $T_{l(0)}M$.
\end{proposition}
{\it Proof: }
For notational convenience, we denote 
$$\tilde l^+ = \tilde l, \tilde l^- = \tilde l', H^+ = H, H^- = H',  \phi^+_t = \phi_t \text{ and } \phi^-_t = c \circ \phi_{1-t} \circ \phi_1^{-1} \circ c \text{ for } t \in [0,1],$$
then $\phi^\pm$ is the flow generated by $H^\pm$.
Assume that we can choose the trivialization $\Phi_z : T_{v(z)}M \to \CC^n$ of $v^*TM$ so that $\Phi_{\bar z} = c_z \circ \Phi_z \circ dc$, where $c_z : \CC^n \to \CC^n$ is the complex conjugation, which takes $\omega_0$ to $-\omega_0$. In particular, $\Phi_r(T_{v(r)}L) = c_z \circ \Phi_r \circ dc (T_{v(r)}L) = \RR^n$ for $r \in [-1, 1]$. Define the following paths of symplectic matrices:
\begin{equation*}
F_t = \Phi_{e^{i\pi t}} \circ d\phi_t \circ \Phi_1^{-1} \text{ for } t \in [0,2] \text{ and } 
F^{\pm}_t  = \Phi_{\pm e^{i\pi t}} \circ d\phi^{\pm}_t \circ \Phi_{\pm 1}^{-1} \text{ for } t \in [0,1],
\end{equation*}
Then $F_t = c_z\circ F_{2-t} \circ F_2^{-1} \circ c_z$ and 
$$\mu_{\HH}(\tilde \gamma) = \mu((F_t, \1) \triangle, \triangle) \text{ and } \mu_{H^\pm}(\tilde l^\pm) = \mu(F^\pm_t\RR^n \dsum \RR^n, \triangle)$$
where $\triangle : \CC^n \to \CC^n \dsum \CC^n$ is the diagonal and the symplectic structure on $\CC^n \dsum \CC^n$ is given by $\Omega_0 = \omega_0 \dsum (-\omega_0)$. We have by additivity of Maslov index:
$$\mu_{\HH}(\tilde \gamma) = \mu((F^+_t, \1)\triangle, \triangle) + \mu((F^-_t \circ F_1, \1)\triangle, \triangle)$$
and the left hand side of \eqref{involution:indexdiffeq} is the sum of the following differences:
$$\mu(F^+_t\RR^n \dsum \RR^n, \triangle) - \mu((F^+_t, \1)\triangle, \triangle) \text{ and } \mu(F^-_t\RR^n \dsum \RR^n, \triangle) - \mu((F^-_t \circ F_1, \1)\triangle, \triangle).$$

For $F \in Sp(\CC^{n})$, $(F, \1)^{-1}\triangle = (\1, F)\triangle$, thus the first difference is
\begin{equation*}
\begin{split}
& \mu(F^+_t\RR^n \dsum \RR^n, \triangle) - \mu((F^+_t, \1)\triangle, \triangle) = \mu((\1, F^+_t)\triangle, \triangle) - \mu((\1, F^+_t)\triangle, \RR^n\dsum \RR^n)\\
= & s(\RR^n\dsum \RR^n, \triangle; \triangle, (\1, F_1)\triangle) = s(\RR^n\dsum \RR^n, (\1, F_1)\triangle; \triangle, (\1, F_1)\triangle)
\end{split}
\end{equation*}
where $s$ is the H\"omander index 
(cf. \cite{RobbinSalamon2}) and the last equality follows from the following properties of H\"ormander index for Lagrangian subspaces $A, B, C, D, D'$:
\begin{equation}\label{invol:hormander}
\begin{split}
s(A, B; A, C) = s(A, B; A, C) - & s(A, C; A, C) = s(C, B; A, C) \\
s(A, B; C, D) - s(A, B;& C, D') = s(A, B; D', D).
\end{split}
\end{equation}
Let $c': \CC^n \dsum \CC^n \to \CC^n \dsum \CC^n : (z_1, z_2) \mapsto (z_2, z_1)$, then $c'$ preserves $\triangle$ and $\RR^n \dsum \RR^n$ while reverses the sign of the symplectic structure, thus
$$s(\RR^n \dsum \RR^n, (\1, F_1)\triangle; \triangle, (\1, F_1)\triangle) = - s(\RR^n \dsum \RR^n, (F_1, \1)\triangle; \triangle, (F_1, \1)\triangle).$$
For the second difference, 
we get 
\begin{equation*}
\begin{split}
& \mu(F^-_t\RR^n \dsum \RR^n, \triangle) - \mu((F^-_t \circ F_1, \1)\triangle, \triangle) \\
= & \mu((\1, F^-_t)\triangle, (F_1, \1)\triangle) - \mu((\1, F^-_t)\triangle, \RR^n \dsum \RR^n)\\
= & s(\RR^n \dsum \RR^n, (F_1, \1) \triangle; \triangle, (\1, F^-_1)\triangle).
\end{split}
\end{equation*}

It follows that the difference 
on the left side of \eqref{involution:indexdiffeq} is
\begin{equation*}\begin{split}
& s(\RR^n \dsum \RR^n, (F_1, \1)\triangle; \triangle, (\1, F^-_1)\triangle) - s(\RR^n \dsum \RR^n, (F_1, \1)\triangle; \triangle, (F_1, \1)\triangle) \\
= & s(\RR^n \dsum \RR^n, (F_1, \1)\triangle; (F_1, \1)\triangle, (\1, F^-_1)\triangle).
\end{split}
\end{equation*}

We now identify the last H\"omander index as the signature.
Let $L = \RR^n \dsum \RR^n$, $K = (F_1, \1)\triangle$ and $L' = (\1, F^-_1)\triangle$, then they are pairwisely transverse, by the non-degeneracy assumption. Thus in the splitting $\CC^{2n} = L \dsum K$ we may write $L'$ as the graph of an invertible linear map $f : K \to K^* \simeq L$ and let $\bar{KL'} = graph(tf)$, $t \in [0,1]$ be the path of Lagrangian subspaces connecting $K$ to $L'$ then
\begin{equation}\label{invol:indproof}
s(L, K; K, L') = \mu(\bar{KL'}, K) - \mu(\bar{KL'}, L) = \mu(\bar{KL'}, K) = \frac{1}{2} sign(Q')
\end{equation}
where $Q'(v) = \Omega_0(v, f(v))$ for $v \in K$ is a quadratic form on $K$. Choose the following coordinates
\begin{equation*}\begin{split}
& L = \{(x, y) | x, y \in \RR^n\}, K = \{(F_1(z), z) | z \in \CC^n\}\text{ and } \\
& L' = \{(\bar w, F^-_1(\bar w)) | w \in \CC^n\} = \{(\bar w, \bar{F_1^{-1}(w)})\} = \{(\bar{F_1(w)}, \bar w)\},
\end{split}
\end{equation*}
where we note $F^-_1 = c\circ F_1^{-1} \circ c$, then it's easy to check that 
$$f : K \to L : z \mapsto (x, y) = -(F_1(z) + \bar{F_1(z)}, z + \bar z)$$
and for $v = (F_1(z), z)$
\begin{equation*}\begin{split}
Q'(v) = & -\omega_0(F_1(z), \bar{F_1(z)}) + \omega_0(z, \bar z) \\
= & -\omega(F_1(z), F_1\circ F_2^{-1}(\bar z)) + \omega_0(z, \bar z) \\
= & \omega_0((\1 - F_2)(z) , \bar z) \\
= & Q(z).
\end{split}
\end{equation*}
Together with \eqref{invol:indproof}, we are done.

We now show the existence of a trivialization $\Phi_z$ with $\Phi_{\bar z} = c_z \circ \Phi_z \circ dc$. Let $V^\pm$ be the $\pm 1$ eigen-bundle of $dc$ action on $v|_{[-1, 1]}^*TM$, then they are transversal Lagrangian subbundles. Since $[-1, 1]$ is contractible, we trivialize $V^+$ and choose a section $\{e^j_r\}_{j = 1}^n$ for $r \in [-1, 1]$ of the frame bundle. The induced trivialization of $V^-$ is then given by $\{f^j_r\}_{j = 1}^n$ where $\omega(e^j_r, f^k_r) = \delta_{kj}$. Then the trivialization $\Phi_r$ can be defined by $\{e^j_r, 
f^k_r\} \mapsto$ standard basis of $\CC^n = \RR^n \dsum i \RR^n$. 
Then the trivialization $\Phi_r$
satisfies $\Phi_r = c_z \circ \Phi_r \circ dc$. Extend it to $D^2_+$ to obtain trivialization $\Phi_z$ for $z \in D^2_+$. Now define $\Phi_z$ for $z \in D^2_-$ by $\Phi_z = c_z \circ \Phi_{\bar z} \circ dc$ and $\Phi_{z \in D^2}$ gives a continuous trivialization of $v^*TM$ with the desired property.
\qed

\subsection{Diagonal}

For $(M, L) = (X \times X, \triangle)$, the doubling construction applies. Let $p_i : M \to X$, for $i = 1, 2$, be the projection to the $i$-th factor, then we obtain the following maps
$$\delta_i = p_i \circ \delta : Map(\Sigma, \partial \Sigma; M, \triangle) \to Map(\Sigma \union_\partial \bar\Sigma; X)$$
which are natural isomorphism between the spaces of continuous maps. 
As special cases, the doubling gives isomorphisms of the path / loop spaces and the respective covering spaces:
$$\delta_i : \PPC_\triangle(M) \to \Omega^{(2)}(X) \text{ and }\delta_i : \tilde \PPC_\triangle(M) \to \widetilde \Omega^{(2)}(X)$$
More explicitly, for example, for $l \in \PPC_\triangle(M)$ we write $l(t) = (l_1(t), l_2(t))$ then
$$(\delta_1(l))(t) = 
\left\{
\begin{matrix}
l_1(t) & \text{ for } t \in [0,1] \\
l_2(2-t) & \text{ for } t \in [1,2]
\end{matrix}
\right.$$
This isomorphism extends to their corresponding normed completions as well.
They  also induce the isomorphisms $\delta_i : \pi_2(M, \triangle) \to \pi_2(X)$. The exact sequence of homotopy groups gives
$$\ldots \to \pi_2(\triangle) \to \pi_2(M) \cong \pi_2(X) \times \pi_2(X) \xto j \pi_2(M, \triangle) \to \ldots$$
Then we have for $\beta \in \pi_2(X)$:
$$\delta_1\circ j(\beta, 0) = \delta_2\circ j(0, -\beta) = \beta$$
It's straight forward to see 
that for $\beta \in \pi_2(X)$, $\delta_2 \circ \delta_1^{-1}(\beta) = - \beta = \tau(\beta)$.
The isomorphism of homotopy group gives rise the isomorphism $\delta_i : \Gamma_\triangle \cong \Gamma_\omega$ as well as the corresponding Novikov rings.
More precisely, for $a_\beta e^\beta \in \Lambda_\triangle$, we have
$$\delta_1(a_\beta e^\beta) = a_\beta e^{\delta_1(\beta)} \in \Lambda_\omega \text{ and } \delta_2(a_\beta e^\beta) = (-1)^{\frac{1}{2}\mu_\triangle (\beta)} a_\beta e^{\delta_2(\beta)} \in \Lambda_{-\omega}$$
then $\delta_2\circ \delta_1^{-1} : \Lambda_\omega \to \Lambda_{-\omega}$ coincides with the isomorphism induced by reversing the symplectic structure on $(X, \omega)$ (cf. \cite{HuLalonde} \S 4).

Let $\{H_t, J_t\}_{t \in [0,2]}$ be a pair of periodic Hamiltonian functions and compatible almost complex structures on $(X, \omega)$, then
$$(\bb H_t, \bb J_t) = (H_t \dsum H_{2-t}, \mb{J}_t \dsum -\mb{J}_{2-t})$$
is a $c$-symmetric pair on $M = X \times X$, with symplectic form $\Omega = \omega \dsum (-\omega)$. Let $\{\phi_t\}_{t \in [0,2]}$ denote the Hamiltonian isotopy generated by $H_t$ on $X$, then $\{\psi_t = (\phi_t, \phi_{2-t}\circ \phi_2^{-1})\}_{t \in [0,2]}$ is the Hamiltonian isotopy generated by $\HH_t$ on $M$. It follows that $x \in X$ is a non-degenerate fixed point of $\phi_2$ iff $(x, x) \in \triangle$ is a non-degenerate fixed point of $\psi_2$.

Let $(\HH^1, \bb J^1)$ and $(\HH^2, \bb J^2)$ be the two halves of $(\HH, \bb J)$, i.e.
$$(\HH^1, \bb J^1) = (\HH_t, \bb J_t)_{t \in [0,1]} \text{ and } (\HH^2, \bb J^2) = (\HH_{t+1}, \bb J_{t+1})_{t \in [0,1]}$$
Let $\tilde l \in \tilde \PPC_\triangle M$ be a critical point of $a_{\HH^1}$, then Lemma \ref{involution:actionfunct} implies that $\tilde \gamma = \delta(\tilde l) \in \tilde \Omega^{(2)}(M)$ is a critical point of $a_{\HH}$. Let $\tilde \gamma_1 = p_1(\tilde \gamma) = \delta_1(\tilde l) \in \tilde \Omega(X)$, then it is a critical point of $a_H$. Similarly, $\tilde \gamma_2 = p_2(\tilde \gamma)$ is a critical point of $a_{\underline H}$, with $\underline H_t = H_{2-t}$. Furthermore, the non-degeneracy of any one of these critical points implies that all the rest are also non-degenerate.

\begin{lemma} \label{lem:indexcoincide}
Suppose that all critical points involved are non-degenerate, then $\mu_\HH(\tilde \gamma) = 2\mu_\HH(\tilde l)$.
It follows that 
$$\mu_\HH(\tilde l) = \mu_H(\tilde \gamma_1)$$
\end{lemma}

{\it Proof:}
The critical point $\tilde \gamma$ is determined by it projection to the two factors, $\tilde \gamma_1$ and $\tilde \gamma_2$. Notice that $\psi_t = (\phi_t, \phi_{2-t} \circ \phi_2^{-1})$, in \eqref{eq:symplecticpath}, the identification $\Phi$ may chosen such that it respects the decomposition $TM = p_1^*TX \dsum p_2^*TX$. Then it's clear that 
$$\mu_\HH(\tilde \gamma) = \mu_H(\tilde \gamma_1) + \mu_{\underline H}(\tilde \gamma_2)$$
Similar to Lemma 5.2 of \cite{HuLalonde}, straight forward computation shows that
$$\mu_H(\tilde \gamma_1) = \mu_{\underline H}(\tilde \gamma_2) \Rightarrow \mu_\HH(\tilde \gamma) = 2 \mu_H(\tilde \gamma_1)$$

Now we only have to see that $\mu_\HH(\tilde \gamma) = 2\mu_\HH(\tilde l)$. By Lemma \ref{lem:reflectionindex} and Proposition \ref{invol:indexdiff}, we only need to compute $sign(Q)$. Let $\gamma(0) = (x, x) \in \triangle$ and $\xi_1, \xi_2 \in T_x X$, then $\xi = (\xi_1, \xi_2) \in T_{\gamma(0)} M$ and
\begin{equation*}
\begin{split}
& Q(\xi, \xi)  = \Omega((\1 - d\psi_2)(\xi_1, \xi_2), (\xi_2, \xi_1))\\
= & \omega((\1-d\phi_2)\xi_1, \xi_2) - \omega((\1 - d\phi_2)\xi_2, \xi_1) \\
= & 2\omega((\1 - d\phi_2)\xi_1, \xi_2) 
\end{split}
\end{equation*}
It follows that $sign(Q) = 0$.
\qed

\subsection{Proof of the lemma}
The lemma follows from the following proposition and Proposition 4.2 of \cite{HuLalonde} which relates the quantum homology of opposite symplectic structures. 
\begin{proposition}\label{diagonal:regularity}
$\delta_1$ induces a natural isomorphism of the Floer theories 
$$\delta_1 : FH_*(M, \triangle; \Omega) \simeq FH_*(X, \omega).$$
\end{proposition}
\noindent
{\it Proof:}
Using the notations from the last subsection, we first compare the action functionals.
Let $\tilde l = [l, w] \in \tilde\PPC_\triangle M$ and $\tilde \gamma_1 = [\gamma_1, v_1]$ so that $\tilde \gamma_1 = \delta_1(\tilde l)$, then
$$a_H([\gamma_1, v_1]) = - \int_{D^2} v_1^*\omega + \int_{[0,2]} H_t(\gamma_1(t)) dt = - \int_{D^2_+} w^*\Omega + \int_{[0,1]} \HH_t(l(t)) dt = a_{\HH}([l, w]).$$
Let $\{\xi_t\}_{t \in [0,2]}$ be a vector field along $\gamma_1$, then $\{\eta_t = (\xi_t, \xi_{2-t})\}_{t \in [0,1]}$ is a vector field along $l$ with $\eta_{0,1} \in T\triangle$ and vice versa. This gives the isomorphism on the tangent spaces:
$$D\delta_1 : T_l \PPC_\triangle M \to T_{\gamma_1} \Omega^{(2)}(X) : \eta \mapsto \xi.$$
It then follows that for $\eta, \eta' \in T_l \PPC_\triangle M$ and the corresponding $\xi$'s:
$$(\xi, \xi')_{\mb J} = \int_{[0,2]} \omega(\xi_t, J_t(\xi'_t)) dt = \int_{[0,1]} \omega(\xi_t, J_t(\xi'_t)) dt + \omega(\xi_{2-t}, J_{2-t}(\xi'_{2-t})) dt = (\eta, \eta')_{\bb J}.$$
From these we see that the Floer equations for the two theories are identified by $\delta_1$ and the moduli spaces of smooth solutions are isomorphic for the two theories.

By Lemma \ref{lem:indexcoincide}, the gradings of the two theories coincide via $\delta_1$. We consider the orientations. Let's first orient the moduli spaces of holomorphic discs in $(M, \triangle)$. Here we may assume that the almost complex structures involved are generic. The map $\delta_1$ induces
$$\delta_{1} : H_*(M, \triangle) \to H_*(X)$$
as well as the maps between the moduli spaces of (parametrized) holomorphic objects (discs or spheres):
$$\delta_1 : \widetilde\MMC(M, \triangle; \JJ, B) \to \widetilde \MMC(X; J, \delta_1(B)).$$
The map $\delta_1$ is an isomorphism. We the put the induced orientation on the moduli space of discs. The moduli spaces of caps are similarly related by $\delta_1$ and the orientations for a preferred basis on either theory can be chosen to be compatible with respect to $\delta_1$. It then follows that the orientations of the theories coincide under $\delta_1$.

To identify the two theories in full, we study the compactifications of the moduli spaces, in particular the compactifications by bubbling off holomorphic discs/spheres. The partial compactification given by the broken trajectories is naturally identified by $\delta_1$ and the identification of the Floer equations.

Consider next the moduli spaces of holomorpic discs in $(M, \triangle)$. 
The map $\delta_1$ defined for the moduli spaces above extends to objects with marked points, which, for spheres, are along $\RP^1 \subset \CP^1$ while for the discs, are along the boundary:
$$\delta_1 : \widetilde\MMC_{k}(M, \triangle; \JJ_i, B) \to \widetilde \MMC_{k}(X; J_i, \delta_1(B)) \text{ for } i = 0, 1.$$
When we pass to the unparametrized moduli spaces, we also denote the induced map $\delta_1$.
Next, we consider the evaluation maps from the moduli spaces of objects with $1$-marked point:
$$ev^\triangle: \MMC_1(M, \triangle; \JJ_i, B) \to \triangle \text{ and } ev: \MMC_1(X; J_i, \delta_1(B)) \to X.$$
Let $p_1 : \triangle \to X$ be the natural projection, then we see that 
$$p_1 \circ ev^\triangle = ev \circ \delta_1.$$
In particular, the image of the evaluation map $ev^\triangle$ has at most the same dimension as that of $ev$ (in fact, they coincide via $p_1$):
$$\dim_\RR = 2c_1(TX)(B) + 2n - 4.$$
The bubbling off of spheres are similar. The Floer theory $FH_*(X, \omega)$ is well defined and it follows that $FH_*(M, \triangle; \Omega)$ is well defined as well and they are isomorphic.
\qed

Recall from \cite{HuLalonde} (Proposition 5.5) that the Lagrangian Floer theories of $(M, \triangle, \Omega)$ and $(M, \triangle, -\omega)$ are related by an isomorphism
$$\tau_* : FH_*(M, \triangle, \Omega; \HH, \JJ) \to FH_*(M, \triangle; -\Omega; \underline \HH, \underline \JJ)$$
where $\underline \HH_t = \HH_{2-t}$ and $\underline \JJ_t = -\JJ_{2-t}$ here.
We observe that the involution $c$ on $M$ identifies the tuples:
$$c: (M, \triangle, -\Omega; \underline\HH, \underline\JJ) \to (M, \triangle; \Omega; \HH, \JJ)$$
and the induced map of $c$ on Floer homology composed with $\tau_*$ is the identity map.
Now Lemma \ref{lem:FHrelation} is given by the following diagram
\begin{equation*}
\text{
\xymatrix{
FH_*(M, \triangle; \Omega) \ar[rr]^{\tau_*}_{\simeq} \ar[d]_{\delta_1} && FH_*(M, \triangle; -\Omega) \ar[rr]^{c_*} \ar[d]_{\delta_1} && FH_*(M, \triangle; \Omega) \ar[lld]^{\delta_2}\\
FH_*(X,\omega) \ar[rr]^{\tau}_{\simeq} && FH_*(X, -\omega)
}
}
\end{equation*}
The commutativity of the left square follows from the discussion of reversing the symplectic structure in \cite{HuLalonde} (\S 4 -- 5), while it's obvious that the right triangle commutes.

\begin{corollary}
The half pair of pants product is well defined for $FH_*(M, \triangle)$ and it has a unit.
\end{corollary}
{\it Proof:}
Everything is induced from $FH_*(X,\omega)$ using the map $\delta_1$.
\qed

\section{Seidel elements and the Albers map}\label{elements}
Let $\Omega_0\Ham(M, \Omega)$ be the space of loops in $\Ham(M, \Omega)$ based at $\1$. It's a group under pointwise composition. 
In $\Omega_0\Ham(M, \Omega)$, a loop $g$ is \emph{split} if $g = (g_1, g_2)$ is in the image of the natural maps
$$\Omega_0\Ham(X, \omega) \times \Omega_0 \Ham(X, -\omega) \to \Omega_0\Ham(M, \Omega)$$
Otherwise, it is \emph{non-split}. Similarly, such notions are defined for the $\pi_1$ of the Hamiltonian groups.
\subsection{Split loops}
In Seidel \cite{Seidel}, the covering space $\tilde \Omega_0\Ham(M, \Omega)$ is defined as
$$\tilde \Omega_0\Ham(M, \Omega) : = \left\{\left.(g, \tilde g) \in \Omega_0\Ham(M, \Omega) \times \Homeo(\tilde \Omega(M)) \right| \tilde g \text{ lifts the action of } g \right\}$$
with covering group $\Gamma_\Omega$. We use $\tilde g$ to denote an element in $\tilde \Omega_0\Ham(M, \Omega)$. The results in \cite{HuLalonde} imply that, similar to \cite{Seidel}, $\tilde g$ defines a homomorphism $FH_*(\tilde g)$ of $FH_*(M, \triangle)$ as a module over itself. Recall that $\delta_1: \Gamma_\triangle \cong \Gamma_\omega$. Moreover, in the homotopy exact sequence
$$\ldots \to \pi_2(\triangle) \xto i \pi_2(M) \to \pi_2(M, \triangle) \to \ldots$$
we have $\img(i) \subset \ker c_1 \inter \ker \Omega$, from which it follows that $\Gamma_\Omega \cong \Gamma_\triangle$.

In the following, we parametrize the loops in $\Omega_0\Ham(X, \omega)$ by $[0,2]$ and those in $\Omega_0\Ham(M, \Omega)$ by $[0,1]$. For $\alpha \in \Omega_0\Ham(X, \omega)$, define the reparametrization $\alpha^{(\frac{1}{2})}(t) = \alpha(2t)$ for $t \in [0,1]$. The natural injective map
$$i_+: \Omega_0\Ham(X, \omega) \to \Omega_0\Ham(M, \Omega) : \alpha \mapsto \alpha_+ = (\alpha^{(\frac{1}{2})}, \1)$$
lifts to an injective map $\tilde i_+$ on the corresponding covering spaces (see the proof of Proposition \ref{prop:deltacommute}). For $\tilde \alpha \in \tilde \Omega_0\Ham(X, \omega)$, let $\tilde \alpha_+ = \tilde i_+(\tilde \alpha_+) \in \tilde \Omega_0\Ham(M, \Omega)$ and $\tilde \alpha_- = \tilde i_-(\tilde \alpha)$ where $\tilde i_-$ is the lifting of
$$i_-: \Omega_0\Ham(X, -\omega) \to \Omega_0\Ham(M, \Omega) : \alpha \mapsto \alpha_- = (\1, (\alpha^-)^{(\frac{1}{2})})$$
We note that $\tilde \alpha_\bullet$ is determined by the image of any element in $\tilde \Omega_0(M)$ by the unique lifting property of covering space. Take the trivial loop $p = (x, y) \in M$, then $x \in M$ is a trivial loop in $\Omega_0(X)$. Let $\tilde \alpha(\tilde x) = [\alpha(x), w] \in \tilde \Omega_0(X)$, where $\tilde x = [x, x] \in \tilde \Omega_0(X)$. Then $\tilde \alpha_+(\tilde p) = [(\alpha^{(\frac{1}{2})}(x), y), w \times \{y\}]$

\begin{proposition}\label{prop:deltacommute}
The following diagram commutes
\begin{equation*}
\text{
\xymatrix{
FH_*(M, \triangle; \Omega) \ar[rr]^{\delta_1} \ar[d]_{FH_*(\tilde \alpha_\pm)} && FH_*(X, \omega) \ar[d]^{FH_*(\tilde \alpha)} \\
FH_*(M, \triangle; \Omega) \ar[rr]^{\delta_1} && FH_*(X, \omega) \\
}
}
\end{equation*}
A similar diagram is commutative with $\delta_2$ and $FH_*(\tilde \alpha^-)$ in places of $\delta_1$ and $FH_*(\tilde \alpha)$.
\end{proposition}
{\it Proof:}
We describe the case for $\tilde \alpha_+$ and $\tilde \alpha_-$ is similar.
Let $\tilde l \in \tilde \PPC_\triangle M$ and $\tilde \gamma = \delta_1(\tilde l) \in \tilde \Omega^{(2)}(X)$.
By definition we have $l(t) = (\gamma(t), \gamma(2-t))$ for $t \in [0,1]$ and $h_1$ acts on $l$ by
$$(\alpha_+ \circ l)(t) = (\alpha_{2t} \circ \gamma(t), \gamma(2-t))$$
Then 
$$(\delta_1( \alpha_+ \circ l))(t) = 
\left\{\begin{matrix}
       \alpha_{2t}(\gamma(t)) & \text{ for } t \in [0,1] \\
       \gamma(t) & \text{ for } t \in [1,2]
       \end{matrix}
\right.$$
which implies that
$$\delta_1(\alpha_+\circ l) = (\alpha^{(\frac{1}{2})} \# \1)\circ \gamma = (\alpha^{(\frac{1}{2})}\# \1) \circ \delta_1 (l)$$
Notice that $\alpha^{(\frac{1}{2})}\# \1$ and $\alpha$ differ by a reparametrization. The equality above lifts to the covering of the loop spaces and gives a chain level identity for the respective Floer theories. In particular
$$\delta_1 \circ FH_*(\tilde \alpha_+) = FH_*(\tilde \alpha) \circ \delta_1$$
\qed

For $\alpha \in \Omega_0\Ham(X, \omega)$, let $\tilde \alpha$ be a lifting to $\tilde \Omega_0\Ham(X, \omega)$. The corresponding Seidel element is
$$\Psi^X(\tilde \alpha) := FH_*(\tilde \alpha)(\1) \in FH_*(X, \omega)$$
where $\1$ is the unit of the pair of pants product.
Moreover, for any other lifting $\tilde \alpha'$ of $\alpha$, there is $B \in \Gamma_\omega$ such that 
$$\Psi^X(\tilde \alpha') = e^B \Psi^X(\tilde \alpha)$$
Similarly, the Lagragian Seidel element is given by
$$\Psi^\triangle(\tilde \alpha_+) = \Psi^\triangle(\tilde \alpha_-) = FH_*(\tilde \alpha_+)(\1)$$
where $\1$ is the unit of the half pair of pants product.
\begin{corollary}\label{coro:splitSeidel}
For $\tilde \alpha$, $\tilde \alpha_\pm$ as given above we have $\delta_1(\Psi^\triangle(\tilde \alpha_+)) = \delta_1(\Psi^\triangle(\tilde \alpha_-)) = \Psi^X(\tilde \alpha)$. \qed
\end{corollary}
\noindent
Since any split loop is the product of $(\alpha, \1)$ and $(\1, \alpha')$, it follows that the split loops in $\Omega_0\Ham(M, \Omega)$ generate Seidel elements in $FH_*(X, \omega)$.

\subsection{The Albers' map}
Here we argue that the Albers' map is well defined for the example under consideration, where $X = (S^2 \times S^2, \omega_0 \dsum \lambda \omega_0)$ with $\lambda \in (1, 2]$. Recall that the map $\ms A : FH_*(M, \Omega) \to FH_*(M, \triangle; \Omega)$ is defined by counting of maps from the ``chimney domain'' $\RR \times [0,1] / \sim$:
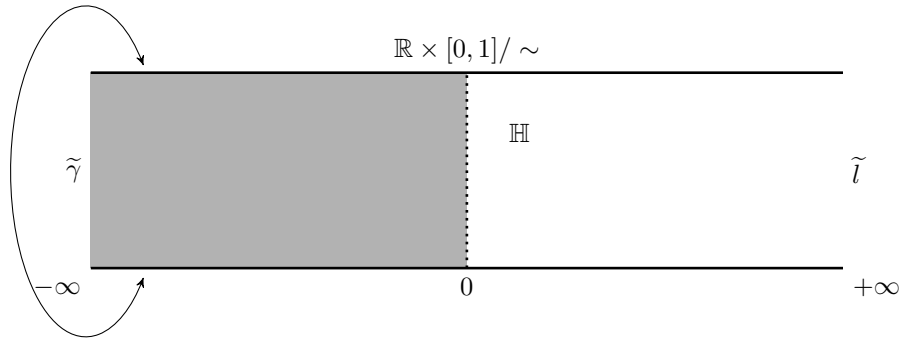
\begin{center}
\begin{tikzpicture}[scale=1,>=stealth']
\fill[color=gray!60!white] (-5,-1.3) -- (0,-1.3) -- (0,1.3) -- (-5,1.3) -- cycle;
\draw[line width=1pt] (-5,-1.3) node[below left,name=ninf] {$-\infty$} -- (0,-1.3) node[below, name=origin] {$0$} -- (5,-1.3) node[below right, name=pinf] {$+\infty$};
\draw[line width=1pt] (-5,1.3) node[name=ninfup] {} -- (0,1.3) node[name=originup] {} -- (5,1.3) node[name=pinfup] {};
\draw[dotted, line width = 1pt] (0,-1.3) -- (0,1.3);
\draw (0.7,0.5) node {$\HH$};
\draw (-5,0) node[left] {$\tilde \gamma$}
      (5,0) node[right] {$\tilde l$}
      (0,1.6) node {$\RR \times [0,1] / \sim$};
\draw[<->] (-4.3,1.4) arc (40:320:1 and 2.2);
\end{tikzpicture}
\end{center}
where $(s, 0) \sim (s, 1)$ when $s \varleq 0$, and the conformal structure at $(0, 0)$ is given by $\sqrt{z}$. In the figure above, the shaded left half of the strip has its two boundaries glued together forming a half infinite cylinder. At $-\infty$ it converges to $\tilde \gamma$, a critical point for the Floer theory $FH_*(M, \Omega)$, while at $+\infty$ it converges to $\tilde l$, a critical point for the Floer theory $FH_*(M, \triangle; \Omega)$.

In \cite{Albers}, the map $\ms A$ is defined for monotone Lagrangians. Here, $(M, \triangle)$ is not monotone because
$$c_1(TM)((01\bar{00})-(10\bar{00})) = 0 \text{ while } \omega((01\bar{00)}-(10\bar{00})) = \lambda - 1 > 0$$
On the other hand, for generic $\omega$-compatible $J$ on $X$, the class $(01\bar{00})-(10\bar{00})$ is not represented by $J$-holomorphic spheres. In fact, the space of non-generic $J$'s has codimension $2$. We choose such a generic pair $(\mb H, \mb J)$ (for the Floer theory $FH_*(X, \omega)$) then the corresponding $c$-symmetric pair $(\HH, \JJ)$ on $(M, \Omega)$ is also generic for the Floer theories $FH_*(M, \Omega)$ and $FH_*(M, \triangle; \Omega)$. 
Since there is no holomorphic disc with non-positive Maslov number, the compactification of the $0$-dimensional ``chimney'' moduli spaces would not contain disc bubblings. Similarly, we see that sphere bubblings can also be ruled out. It then follows that the map $\ms A$ is well-defined.

\subsection{Non-split loops}

We showed that in $\Omega_0\Ham(M, \Omega)$, there could be non-split loops, by computing directly the corresponding Seidel elements in $QH_*(M, \Omega)$. For such loops, Proposition \ref{prop:deltacommute} does not apply. On the other hand, let $g \in \Omega_0\Ham(M, \Omega)$ be a non-split loop and $\tilde g$ be a lifting to $\tilde \Omega_0\Ham(M, \Omega)$, then it defines a Seidel element $\Psi^M(\tilde g) \in FH_*(M, \Omega)$. The Albers' map \cite{Albers} $\ms A$ relates $FH_*(M, \Omega)$ and $FH_*(M, \triangle)$ when it's well defined, in which case, we have
$$\ms A \circ \Psi^M(\tilde g) = \Psi^\triangle(\tilde g) \in FH_*(M, \triangle) \text{ and } \delta_1\circ \ms A \circ \Psi^M(\tilde g) \in FH_*(X, \omega)$$

Consider $(X, \omega) = (S^2 \times S^2, \omega_0 \dsum \lambda \omega_0)$ for $\lambda \in (1, 2]$ and compute $\delta_1\circ \ms A \circ \Psi^M(\tilde g)$ for a non-split loop $g$. 
Also recall that 
McDuff \cite{McDuff} showed that liftings of the loops in the Hamiltonian group may be chosen such that
$$\Psi^X : \pi_1\Ham(X, \omega) \to QH_*(X, \omega) : \alpha \mapsto \tilde \alpha \mapsto \Psi^X_\alpha := \Psi^X(\tilde \alpha)$$
is a group homomorphism. 
Let $\psi = \Psi^M_{ S'}$, i.e.
$$\psi = \left[(01\bar 1 \bar 1)-(11\bar 1 \bar 0)\right]e^{\frac{1}{2}(1000) + h\left[(0001) + (1000)\right]}$$
To compute $\delta_1\circ \ms A(\psi)$, we note first that $\delta_1\circ \ms A$ is linear with respect to the identifications of the Novikov rings. 
Consider the following Seidel elements of split loops:
$$\Psi^M_{R_1} = (01\bar 1 \bar 1) e^{\frac{1}{2}(1000)} \text{ and } \Psi^M_{\bar R_2} = -(11\bar 1\bar 0) e^{-\frac{1}{2}(0001)}$$
then 
$$\delta_1 \circ \ms A \circ \Psi^M_{R_1} = \delta_1(\Psi^\triangle(R_1)) = \Psi^X(r_1) = (01)e^{\frac{1}{2}(10)}$$
$$\delta_1 \circ \ms A \circ \Psi^M_{\bar R_2} = \delta_1(\Psi^\triangle(\bar R_2)) = \Psi^X(r_2) = (10)e^{\frac{1}{2}(01)}$$
where we use $r_i$ to denote the rotation of the $i$-th $S^2$ factor of $X$. In particular, we recall that via the identifications of Novikov rings,
$$e^{\frac{1}{2}(1000)} \mapsto e^{\frac{1}{2}(10)} \text{ and } e^{-\frac{1}{2}(0001)} \mapsto e^{\frac{1}{2}(10)}$$
It follows that
$$\delta_1\circ\ms A(\psi) = \left[(01) + (10)\right] e^{\frac{1}{2}(10) + h\left[(10)-(01)\right]} = \Psi^X(s)$$
where $s$ represents the element of infinite order in $\pi_1\Ham(X, \omega)$.
In summary, we showed
\begin{proposition}
Under $\delta_1\circ \ms A$ the Seidel elements of the non-split loops in $\Ham(M, \Omega)$ constructed in \cite{HuLalonde3} map to the Seidel elements of the loops of infinite order in $\Ham(X, \omega)$.
\qed
\end{proposition}

\label{lastpage-01}
\end{document}